\documentclass[12pt]{amsart}
\usepackage{amssymb}
\usepackage{amsthm}
\usepackage{latexsym}
\usepackage{amsmath}
\usepackage{amscd}
\usepackage{bm}
\usepackage{eufrak}
\usepackage[all]{xy}
\usepackage[usenames]{color}
\usepackage[dvips]{graphicx}

\theoremstyle{plain}
\newtheorem{teo}{Theorem}[section]

\newtheorem{prop}[teo]{Proposition}
\theoremstyle{definition}

\newtheorem{obs}{Remark}[section]

\DeclareMathAlphabet{\mathpzc}{OT1}{pzc}{m}{it} \DeclareMathOperator{\fix}{{Fix}} \DeclareMathOperator{\aut}{{Aut}}

\newcommand{\R}{\mathbb R}

\def\Z{{\mathbb Z}}
\def\H{{\mathbb H}}
\def\F{{\mathbb F}}

\def\P{{\mathbb P}}
\def\C{{\mathbb C}}
\def\Z{{\mathbb Z}}
\begin{document}

\bibliographystyle{amsplain}

\title[Construcci\'on]{Riemann Surfaces of genus $g$ with an automorphism of order $p$ prime and $p>g$}
\author{\textit{Giancarlo Urz\'ua}}

\email{gian@umich.edu}

\maketitle

\begin{abstract}
The present work completes the classification of the compact Riemann surfaces of genus $g$ with an analytic automorphism of order $p$ (prime
number) and $p>g$. More precisely, we construct a parameterization space for them, we compute their groups of uniformization and we compute
their full automorphism groups. Also, we give affine equations in $\C^2$ for special cases and some implications on the components of the
singular locus of the moduli space of smooth curves of genus $g$.
\end{abstract}

\section{Introduction}\label{s1}

We first fix the notations.\\
$S$ compact Riemann surface or smooth projective curve over $\C$\\
$p$ prime number, $g$ genus of $S$\\
$\aut(S)= \{ f : S \longrightarrow S: f $ analytic automorphism$\}$ (full automorphism group of $S$)\\
$\Z/n\Z$ cyclic group of order $n$\\
$D_m$ the dihedral group of order $2m$\\
For a set $A$, $\vert A \vert$ is its cardinality\\
If a group $G$ acts on a set, $\fix(G)$ is the set of the fixed points by the action\\
$\Delta = \{ z \in \C : \vert z \vert < 1\}$ the unit disk\\
$\langle T \rangle$ is the subgroup generated by $T$\\
$\P^1$ Riemann sphere or one dimensional projective space over $\C$\\
$S \simeq S'$ means that $S$ and $S'$ are isomorphic as Riemann surfaces\\
$M_g$ denotes the moduli space of smooth curves of genus $g$\\

We begin with the following theorem (see \cite{F-K} chapter V). If $T \in \aut(S)$ is of prime order $p>g$, then there are only four
possibilities:

\begin{itemize}
\item[1)] $S/<T>$ has genus $1$, $g=1$ and $|\fix(T)|=0$. \item[2)] $S/<T> \simeq \P^1$, $g=0$ and $|\fix(T)|=2$. \item[3)] $S/<T> \simeq \P^1$,
$p=2g+1$ and $|\fix(T)|=3$. \item[4)] $S/<T> \simeq \P^1$, $p=g+1$ and $|\fix(T)|=4$.
\end{itemize}

The last two cases are the non trivial ones. In a more general setting, Solomon Lefschetz studied these surfaces using algebraic equations in
\cite{L}. Among other things, he computed the full automorphism group for surfaces of the third case, which are known as Lefschetz surfaces (see
\cite{Ri-R-1}). Since then, people have worked on Lefschetz surfaces, classifying them completely (for example, see \cite{G} and \cite{Ri-R-1}).
In \cite{G}, the author proves that the Lefschetz surfaces which have as full automorphism group $\Z/p\Z$ are exactly the isolated points of the
singular locus of $M_g$ for $g>3$; he also counted these points for every $g$. Since Lefschetz surfaces are branched on three points on $\P^1$,
their moduli space is zero dimensional. This gives to the problem of classification of Lefschetz surfaces a purely combinatorial character.

For the fourth case, we have one free complex parameter plus the combinatorial problem. In this work, we give a method for working out that
case. As consequences, in section \ref{s8} we construct a parameterization space for them, in section \ref{s9} we compute the full automorphism
group for each of these surfaces, in section \ref{s10} we compute formulas for the number of zero and one dimensional components of the singular
locus of $M_g$ and in section \ref{s11} we write affine equations for families of special cases.

\section{Lefschetz surfaces}\label{s2}

Let $p>3$ be a prime number. We do not consider $p=2$ or $3$ because $p=2$ is not possible, and for $p=3$ we have the unique compact Riemann
surface of genus one admitting an automorphism of order three. We want to show how the combinatorial method, which will be used in a similar way
for the fourth case, can be used to reclassify Lefschetz surfaces. We denote the set of those surfaces by $\bm{\mathpzc{L}_p}$. Let $S
\in{\mathpzc{L}_p}$ and let $T \in{\aut(S)}$ be such that $\langle T \rangle$ has signature $(0;p,p,p)$; this means that $S/\langle T \rangle
\simeq \P^1$ and $|\fix(T)|=3$, so by Riemann-Hurwitz formula $g = \frac{p-1}{2}$.  We refer to \cite{F-K} chapter IV for definitions. Let
$\Gamma$ be a Fuchsian subgroup of $\aut(\Delta)$ such that its signature is $(0;p,p,p)$. This group is unique up to conjugation in
$\aut(\Delta)$. It is known that $\Gamma$ can be presented in the following way.
$$ \Gamma=<x,y \ : \ x^p=y^p=(xy)^p=1 > $$
where $x,y$ are suitable elements in $\aut(\Delta)$. Given this situation, there exists $\Gamma^*$ normal subgroup of $\Gamma$ such that
$\Gamma/\Gamma^* \simeq \Z/p\Z$ and $\Delta/\Gamma^* \simeq S$. Below we show the commutative diagram of coverings induced by these groups.
$$ \xymatrix{ \Delta \ar[d]_{\pi_{\Delta}} \ar[dr]^{\pi} & \\
S \simeq \Delta/ \Gamma^* \ar[r]^{\pi_S} & \Delta/\Gamma \simeq \P^1} $$

\begin{center}
\begin{figure}
\includegraphics[width=13cm]{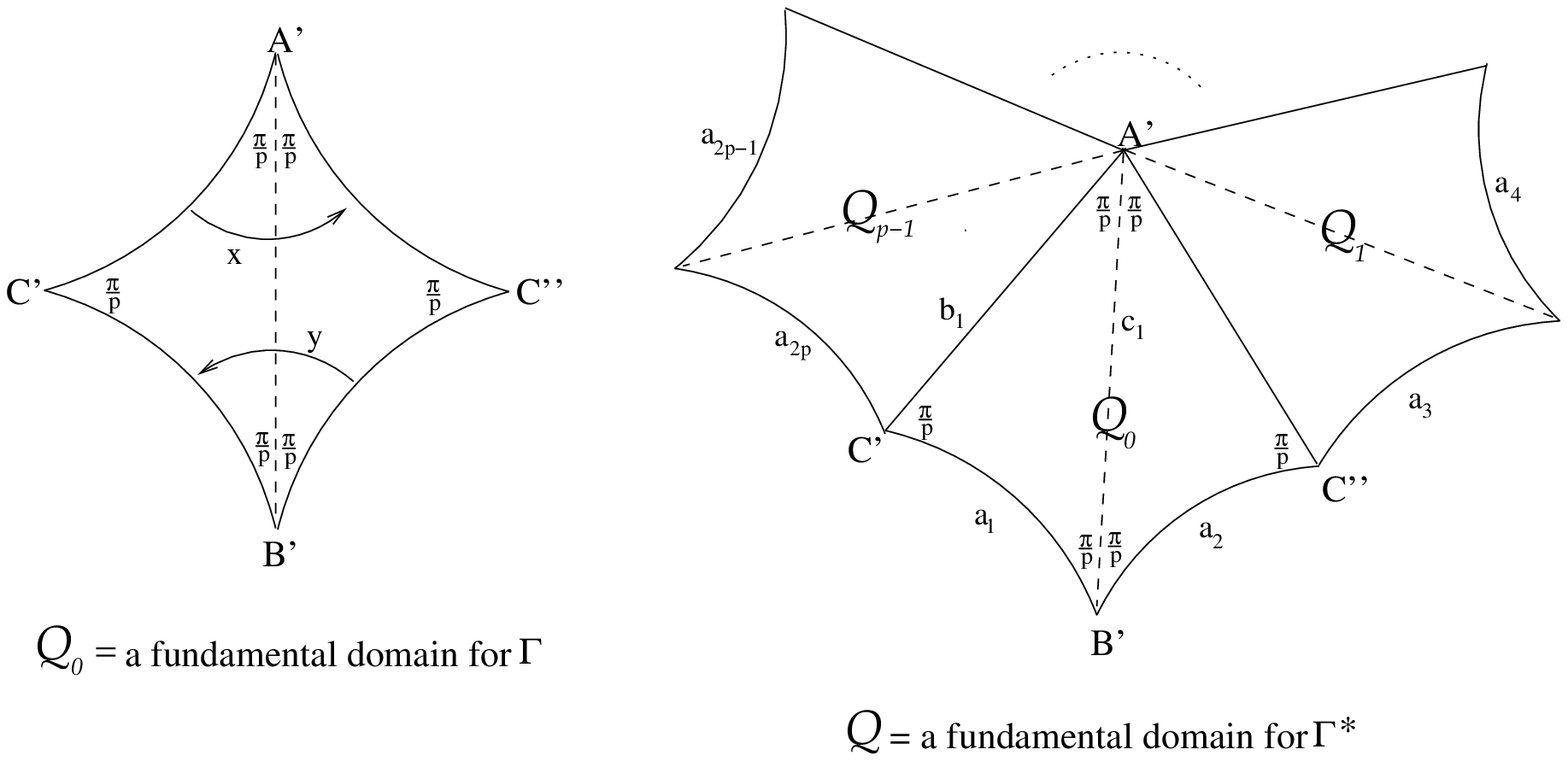}
\caption{Fundamental domains for $\Gamma$ and $\Gamma^*$} \label{f0}
\end{figure}
\end{center}

Let $\{A,B,C \} = \pi_S(\fix(T))$. The following discussion is illustrated in Figure \ref{f0}. We consider as a fundamental domain for $\Gamma$
a quadrilateral $\mathit{Q}_0$ formed by two equilateral triangles with angles $\frac{\pi}{p}$, such that the opposite angles in this
quadrilateral are $\frac{2\pi}{p}$ and $\frac{\pi}{p}$. The vertices corresponding to $\frac{2\pi}{p}$ will be denoted by $A'$ and $B'$, and for
$\frac{\pi}{p}$ we have $C'$ and $C''$, so that $\pi(A')=A$, $\pi(B')=B$ and $\pi(C')=\pi(C'')=C$. We choose to have $x$ as the rotation by
angle $\frac{2\pi}{p}$ and center $A'$, and $y$ as the rotation by angle $\frac{2\pi}{p}$ and center $B'$. For $\Gamma^*$, we consider the
fundamental domain $\mathit{Q}$ given by $p$ copies around $A'$ of $\mathit{Q}_0$. The quadrilaterals $\mathit{Q}_0, \mathit{Q}_1, ...,
\mathit{Q}_{p-1}$ are these copies enumerated counterclockwise and, in the same way, the edges of $\mathit{Q}$ are denoted by $a_1, a_2, ...,
a_{2p}$ such that $a_{2k+1}, a_{2k+2}$ correspond to $\mathit{Q}_k$ for every $k \in \{ 0,1,...,p-1 \}$.  Hence, to obtain the possible groups
$\Gamma^*$ for a fixed $\Gamma$, it is enough to find which are the possible identifications of the edges $a_k$ of $\mathit{Q}$. By the symmetry
of this polygon, it would be enough to show how to glue $a_1$. This is equivalent to looking at all possible surjective homomorphisms $\Gamma
\rightarrow \Z/p\Z \rightarrow 1$ (see \cite{H}).

\begin{prop}
The edge $a_1$ can only be glued with an edge of the form $a_{2k+2}$ for any $k \in \{ 1,2,...,p-2 \}$. If we fix $k$ and $T_{a_1
\leftrightarrow a_{2k+2}}$ is the element in $\aut(\Delta)$ which makes this identification, then $T_{a_1 \leftrightarrow a_{2k+2}} = x^k
y^{-1}$ and $\Gamma^*$ is generated by $\{ x^{jk}  T_{a_1 \leftrightarrow a_{2k+2}} x^{-jk} \}_{j=0}^{p-1}$. \label{p1}
\end{prop}

\begin{proof}
First, we notice that $a_1$ cannot be glued with $a_2$ or $a_{2p}$ because in that case $\Delta / \Gamma^*$ would not be a Riemann surface of
genus $\frac{p-1}{2}$. If we glue $a_1$ with $a_{2k+1}$, then $\Delta / \Gamma$ would not be an oriented Riemann surface. All the other cases
give us Lefschetz surfaces. After that, it is straightforward to compute the generators for $\Gamma^*$.
\end{proof}

\begin{prop}
Let $S,S'$ be in $\mathpzc{L}_p$. If $S \simeq S'$, then there exits an isomorphism $f : S\longrightarrow S'$ such that
$f(\fix(T_S))=\fix(T_{S'})$, where $T_S$ and $T_{S'}$ are automorphisms of order $p$ fixing three points. \label{p2}
\end{prop}

\begin{proof}
First, we want to state Theorem $1$ in \cite{G-G}: Let $X$ be a compact Riemann surface and suppose that $\aut(X)$ contains automorphisms $T_1$,
$T_2$ of same prime order and such that the quotient surfaces $X/\langle T_i \rangle$ ($i=1,2$) are isomorphic to $\P^1$; then $\langle T_1
\rangle$ and $\langle T_2 \rangle$ are conjugate in $\aut(X)$. Let $g: S \longrightarrow S'$ be an isomorphism as Riemann surfaces. Then, $g T_S
g^{-1}$ is an automorphism of $S'$ of order $p$ and so there exists $H \in{\aut(S')}$ such that $Hg T_S (Hg)^{-1}= T_{S'}^{n}$ for some $n
\in{\{1,...,p-1\}}$. Finally, $f=Hg$ satisfies the requirements.
\end{proof}

The last proposition is a key fact in this paper. To know whether two Riemann surfaces $S,S'$ in $\mathpzc{L}_p$ are isomorphic or not, we only
need to observe what is happening around the special points $A$, $B$ and $C$. More precisely, we need to observe how the identification of edges
changes when we look at different fundamental domains $\mathit{Q}$ for a fixed $\Gamma^*$.  By different fundamental domains, we mean domains
which are hyperbolically identical to $\mathit{Q}$, as defined above, but now with center at $B'$ or $C'$ (where those points are preimages of
$B$ or $C$ respectively). Let us denote by $[X,Y]$ a domain identical to $\mathit{Q}$ with center $X'$ such that $\pi(X')=X$ and if $Y'$ is a
vertex opposite to $X'$ (as $B$ was opposite to $A$ for $\mathit{Q}$) , then $\pi(Y')=Y$. If $X,Y \in \{A,B,C \}$, then we have only six
alternatives up to isometries of $\Delta$ for these domains: $[A,B]$, $[A,C]$, $[B,A]$, $[B,C]$, $[C,A]$ and $[C,B]$. Now, according to Prop
\ref{p1}, $\Gamma^*$ is completely determined by a choice of $k$. Hence, after we choose $k$, we will get identifications for all of the domains
$[X,Y]$. Below, we compute $k$ for all of them, $k$ in the sense of Prop \ref{p1}. For a domain $[X,Y]$, the corresponding $k$ will be denoted
by $s[X,Y]$.

\begin{prop}
Given $s[A,B] = k$, then $s[A,C] = p-1-k$,  $s[B,A] = k^{-1}$, $s[B,C] = p-1-k^{-1}$, $s[C,A] = p- (k+1)^{-1}$, $s[C,B] = (k+1)^{-1} - 1$ where
all the numbers involved are taken from $\Z/p\Z$ in $\{1,2,...,p-2 \}$. \label{p3}
\end{prop}

\begin{proof}
Take the triangle of vertices $A'$, $B'$, $C'$ in Figure \ref{f0} whose edges are denoted by $a_1$, $b_1$, $c_1$ as shown. We will first compute
$s[B,C]$. For that, we need to realize how $b_1$ is glued with some $b_{2s[B,C]+2}$. This is equivalent to counting how many commas we have in
the following chain of identifications.
$$ a_1\leftrightarrow a_{2k+2} , a_{2k+1}\leftrightarrow a_{4k+2}, \ldots , a_{2p-2k-1}\leftrightarrow a_{2p} $$
Hence, we obtain the equation $2(s[B,C]+1)k + 2 \equiv 2p (mod \ p)$ and so $s[B,C] \equiv -1-k^{-1} (mod \ p)$. Similarly, for $[C,A]$ we have
$2(s[C,A]+1)s[B,C] + 2 \equiv 2p (mod \ p)$, and so $s[C,A] \equiv - (k+1)^{-1} (mod \ p)$. For the rest, we just apply the operation $\bullet
\mapsto p-1- \bullet$.
\end{proof}

We notice that those numbers are invariants of the conjugacy class of $\Gamma^*$ (and so of the surface it creates), and do not depend on the
choice of points in $\pi^{-1}( \{ A,B,C \})$ which are used to form the domains $[X,Y]$. On the other hand, the set $ \{ 1,2,...,p-2 \}$ is
partitioned  by the set $\{k, p-1-k, k^{-1},  p-1-k^{-1}, p- (k+1)^{-1},  (k+1)^{-1}-1 \}$ (all the numbers taken mod $p$). Therefore, we can
assign to any $S \in{\mathpzc{L}_p}$ the set of numbers
$$ \Omega_k^p = \{k, p-1-k, k^{-1},  p-1-k^{-1}, p- (k+1)^{-1},  (k+1)^{-1}-1 \} $$
for some $k \in{ \{1,2,...,p-2 \}}$ (which is determined by the corresponding $\Gamma^*$).

\begin{teo}
Let $S,S'$ be in $\mathpzc{L}_p$ and let $\Omega_k^p, \Omega_{k'}^p$ be the corresponding sets as above. Then, $S \simeq S'$ if and only if
$\Omega_k^p = \Omega_{k'}^p$. \label{t1}
\end{teo}

\begin{proof}
If $\Omega_k^p = \Omega_{k'}^p$, then it is clear how to build the isomorphism between the surfaces. Now suppose that $S \simeq S'$. Then, by
Prop. \ref{p2}, we know that there exists an isomorphism $f : S\longrightarrow S'$ such that $f(\fix(T_S))=\fix(T_{S'})$, where $T_S$ and
$T_{S'}$ are automorphisms of order $p$ fixing four points. Therefore, both surfaces can be thought as having the same $\Gamma$. Choose special
domains $[X,Y]$ for $S$ and $[X',Y']$ for $S'$ such that $f(X)=X'$ and $f(Y)=Y'$. Then $s[X,Y]=s[X',Y']$ and so, since the sets $\Omega_k^p$
partition $\{ 1,2,...,p-2 \}$, we obtain that the corresponding sets for $S$ and $S'$ are equal.
\end{proof}

\begin{prop}
The following are the only possible cardinalities for $\Omega_k^p$:\\
-$2$, occurring exactly when $p\equiv 1(mod \ 3)$ and $k$ satisfies $k^2+k+1\equiv 0(mod \ p)$.\\
-$3$, occurring exactly for $\Omega_1^p=\{1,\frac{p-1}{2},p-2\} $.\\
-$6$, otherwise. \label{p4}
\end{prop}

\begin{prop}
The set $\{1,2,...,p-2\}$ is partitioned by the sets $\Omega_k^p$ in the following number of subsets: $\frac{p+5}{6}$ if $p\equiv1 (mod \ 3)$,
or $\frac{p+1}{6}$ if $p\equiv-1 (mod \ 3)$. \label{p5}
\end{prop}

The proofs of the last two propositions are elementary and involve only linear equations mod $p$. The following theorem is the result of putting
together Thm. \ref{t1} and Prop. \ref{p5}.

\begin{teo}
The number of analytically distinct compact Riemann surfaces in $\mathpzc{L}_p$ is $\frac{p+5}{6}$ if $p\equiv1 (mod \ 3)$, and $\frac{p+1}{6}$
if $p\equiv-1 (mod \ 3)$. \label{t2}
\end{teo}

A list of the sets $\Omega_k^p$ for small $p$ is given in the appendix. Because of Thm. \ref{t1}, surfaces can be thought as sets $\Omega_k^p$.
As we will see, these sets give enough information to find the full automorphism group and also affine equations for them. It turns out that,
except for the Klein curve $\Omega_2^7$, we have:

\begin{itemize}
\item[-] $\aut(\Omega_k^p) \simeq \Z/p\Z \rtimes \Z/3\Z$ when $\vert \Omega_k^p \vert=2$. \item[-] $\aut(\Omega_1^p) \simeq \Z/2p\Z$. \item[-]
Otherwise, $\aut(\Omega_k^p) \simeq \Z/p\Z$.
\end{itemize}

\begin{proof}
Let $S$ be in $\mathpzc{L}_p$ and $T$ automorphism of order $p$ fixing four points. Consider the following subgroup of $\aut(S)$.
$$\aut(S)' :=\{ R \in{\aut(S)} \ : \ R(\fix(T_S))=\fix(T_S) \ \} $$
One can check that this is the normalizer of $\langle T \rangle$ in $\aut(S)$. Now, we can compute $\aut(S)'$ by using the set $\Omega_k^p$
corresponding to $S$. By keeping track of the domains $[X,Y]$, these sets tell us when we can and cannot identify two domains $[X,Y]$. For
instance, in $\Omega_1^p$ we can identify $[A,B]$ with $[B,A]$ since they have the same $k$, which is $1$. Using this, we can get the above list
for $\aut(S)'$. Now, the question of whether $\langle T \rangle$ is normal or not in $\aut(S)$ is answered by the Singerman's list of finitely
maximal Fuchsian groups (see \cite{S}). This list tell us (in our case) that $\langle T \rangle$ is always normal, except for the case of the
Klein curve.
\end{proof}

An interpretation of this is that the larger the size of $\Omega_k^p$ is, the more rigid the corresponding surface is. Finally, we would like to
explain how to get equations for these surfaces. For this, we will refer to the first part of \cite{G-G1}. Let $C$ be a compact Riemann surface
which has an automorphism $T$ of order $p$ such that $C/ \langle T \rangle \simeq \P^1$. Then, $C$ is isomorphic to the smooth model of the
affine curve given by
$$ y^p = \prod_{i=1}^r (x-a_i)^{m_i} $$
where $a_i$ are distinct complex numbers, $1 \leq m_i \leq p-1$ for every $i$, and $\sum_{i=1}^r m_i \equiv 0 (mod \ p)$.  We take $T$ as $(x,y)
\mapsto (x, e^{\frac{2\pi i}{p}}y)$. The fixed points of $T$ are $\{ (a_1,0),..., (a_r,0) \}$, so the genus $g$ of $C$ is
$\frac{(p-1)(r-2)}{2}$. If $t$ is a local parameter for $C$ in an open set $U$ containing $(a_i,0) $, then $T\vert_U$ is given by $t \mapsto
w^{\sigma_i}t$. The positive integer $\sigma_i$ is usually called rotation number. As a consequence, we get $\sigma_i m_i \equiv 1 (mod \ p)$.
Coming back to Lefschetz surfaces, it is a simple computation that we can choose $\{ 1, k^{-1}, (p-1-k)^{-1} \}$ as rotation numbers around the
three fixed points.  After an automorphism of $\P^1$, we can take  $\{ 0, 1, \infty \}$ as the branch points, obtaining that $ \Omega_k^p $ is
given by $\bm{ y^p=x(x-1)^k}$. For example, the Klein curve $\Omega_2^7$ is $y^7=x(x-1)^2$.

As I said at the beginning, the aim of this section was to illustrate how to reclassify Lefschetz surfaces with our method. Next, we generalize
this method to classify the compact Riemann surfaces which have an automorphism of order $p$ prime fixing four points with quotient $\P^1$.

\section{the theoretical framework for four points}\label{s3}

Now we start the study of the fourth alternative in section \ref{s1}. The set of those surfaces will be denoted by $\bm{\gimel_p}$. The
following theoretical framework is similar to what we had for Lefschetz surfaces. Let $S \in{ \gimel_p }$ and let $T_S$ be an automorphism of
$S$ of order $p$ (prime number) such that $S/  \langle T_S \rangle \simeq \P^1$ and $\vert \fix(T_S) \vert =4$. Notice that, by using the
Riemann-Hurwitz formula, $g=p-1$. We will take $p>2$ because for $p=2$ we have all the compact Riemann surfaces of genus $1$ and their
classification is already well known.  Let $\Gamma$ be a Fuchsian subgroup of $\aut(\Delta)$ such that its signature is $(0;p,p,p,p)$. This
group is unique up to conjugation in $\aut(\Delta)$. It is known that $\Gamma$ can be presented in the following way.
$$ \Gamma=<x,y,z \ : \ x^p=y^p=z^p=(xyz)^p=1 > $$
where $x,y,z$ are suitable elements in $\aut(\Delta)$. Given this situation, there exists $\Gamma^*$ normal subgroup of $\Gamma$ such that
$\Gamma/\Gamma^* \simeq \Z/p\Z$ and $\Delta/\Gamma^* \simeq S$. Below we show the commutative diagram of coverings induced by these groups.
$$ \xymatrix{ \Delta \ar[d]_{\pi_{\Delta}} \ar[dr]^{\pi} & \\
S \simeq \Delta/ \Gamma^* \ar[r]^{\pi_S} & \Delta/\Gamma \simeq \P^1} $$

\begin{center}
\begin{figure}
\includegraphics[width=13cm]{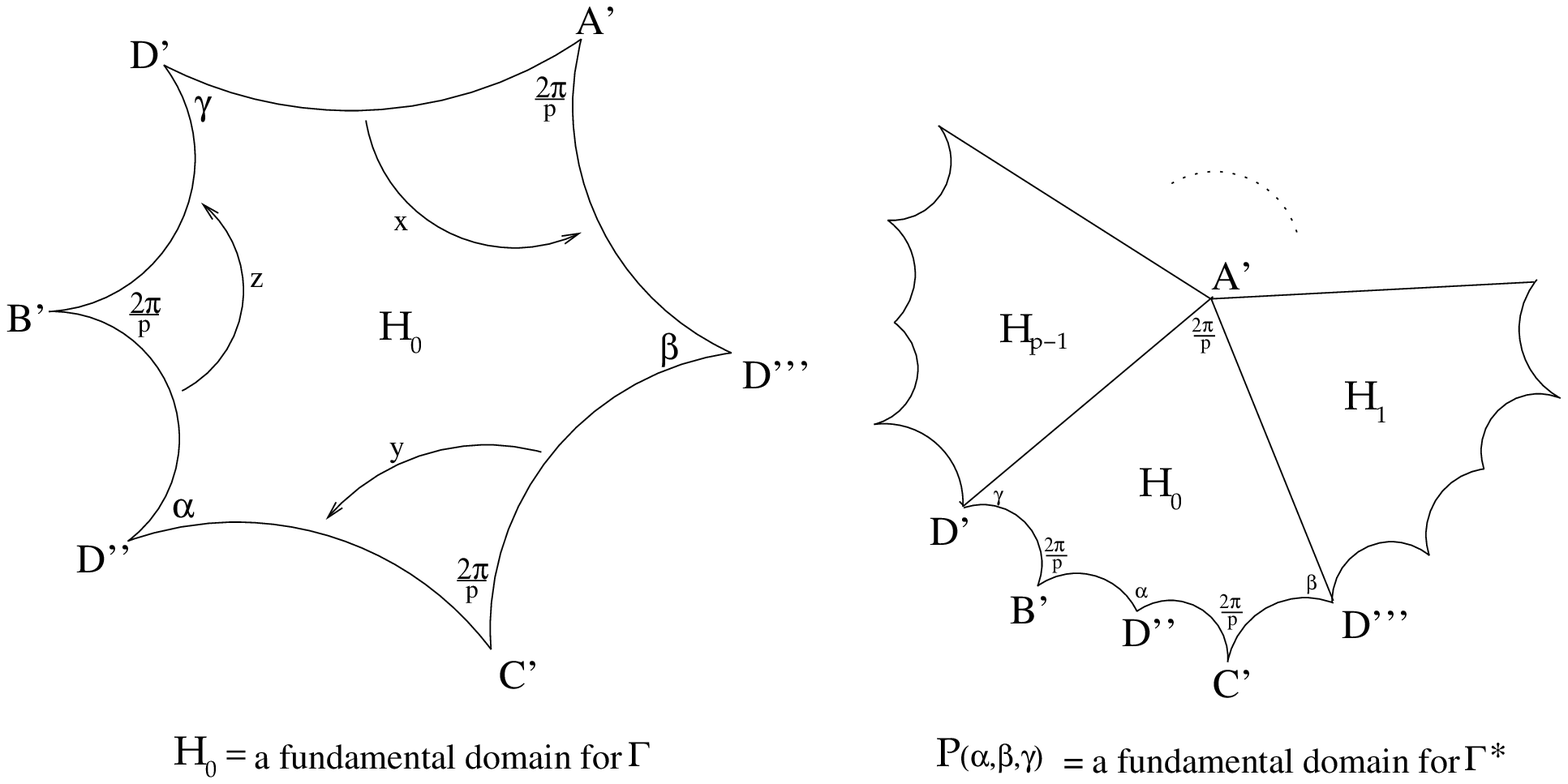}
\caption{Fundamental domains for $\Gamma$ and $\Gamma^*$} \label{f1}
\end{figure}
\end{center}

Let $\{A,B,C,D \} = \pi_S(\fix(T))$. Now, we consider as a fundamental domain for $\Gamma$ the hexagon $\rm{H}_0$ which is shown in Figure
\ref{f1}. The vertices of $\rm{H}_0$ are $\{ A', D', B', D'', C', D''' \}$ where $\pi(A')=A$, $\pi(B')=B$, $\pi(C')= C$ and
$\pi(D')=\pi(D'')=\pi(D''')=D$, and so we have the condition $\alpha + \beta + \gamma = \frac{2 \pi}{p}$. We choose generators of $\Gamma$ $x$,
$y$ and $z$ to be the rotations by angle $\frac{2\pi}{p}$ which fix $A'$, $B'$ and $C'$ respectively (as shown in Figure \ref{f1}). We consider
as a fundamental domain for $\Gamma^*$ the polygon $\rm{P}(\alpha,\beta,\gamma)$ formed by $p$ hexagons identical to $\rm{H}_0$ around $A'$,
denoting them by $\rm{H}_0, \rm{H}_1,...,\rm{H}_{p-1} $. The free complex parameter in $\gimel_p$ is evident in the angular equation $\alpha +
\beta + \gamma = \frac{2 \pi}{p}$.

\section{handling the angular parameter}\label{s4}

\begin{center}
\begin{figure}
\includegraphics[width=8cm]{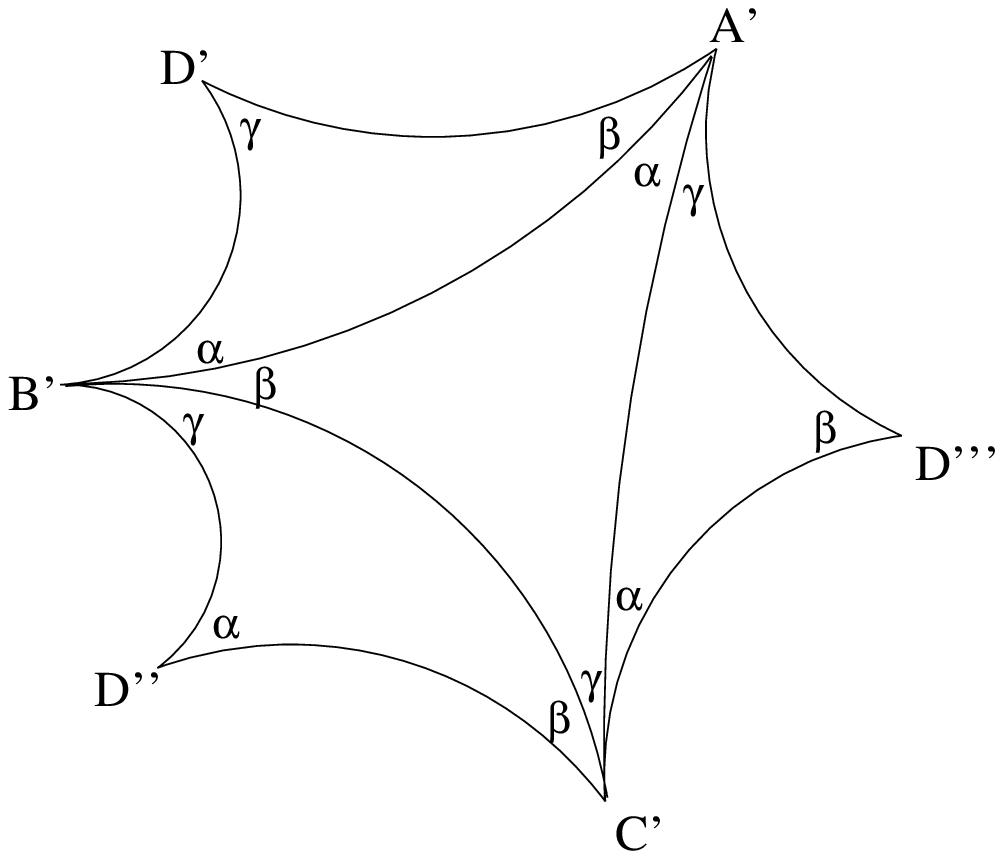}
\caption{$\rm{H}(\alpha,\beta,\gamma)$} \label{f2}
\end{figure}
\end{center}

The purpose of the present section is to manage this free complex parameter. The angular requirement $\alpha + \beta + \gamma = \frac{2 \pi}{p}$
gives the existence of a hyperbolic triangle of angles $\alpha$, $\beta$ and $\gamma$. Then, we can tile $\Delta$ by gluing this triangle on its
edges as shown in Figure \ref{f2}. Consider the hexagon $\rm{H}(\alpha,\beta,\gamma)$ shown in that figure. Since the angular equation above is
satisfied, there are transformations in $\aut(\Delta)$ which give rotations around $A'$, $B'$ and $C'$ in angle $\frac{2 \pi}{p}$, and they
generate a Fuchsian group whose signature is $(0;p,p,p,p)$. Hence, without lost of generality, we can take $\rm{H}_0 = H(\alpha,\beta,\gamma)$.

We denote the tiling on $\Delta$ given by a triangle of angles $\alpha$, $\beta$ and $\gamma$ (in that order) by $\Delta(\alpha,\beta,\gamma)$.
Hence, $\Delta(\alpha,\beta,\gamma) = \Delta(\beta,\gamma,\alpha) = \Delta(\gamma,\alpha,\beta)$. This tiling is unique up to automorphisms of
$\Delta$.  We observe that for a given $\Gamma$ we have infinitely many possible tilings for $\Delta$. Our goal is to find a canonical tiling
$\Delta(\alpha,\beta,\gamma)$ for $\Gamma$ which characterizes this group completely.

\begin{prop}
Let $p$ be an odd prime number and $\Gamma$ be as above. Then, we can find a unique $\Delta(\alpha,\beta,\gamma)$ for $\Gamma$ such that
$0<\alpha \leq \frac{\pi}{p}$ and $0<\beta,\gamma<\frac{\pi}{p}$, except when $\alpha = \frac{\pi}{p}$, for which we have two such tilings (with
the same angles but in opposite order). \label{p6}
\end{prop}

We present a sketch for a proof. First we prove that such a tiling exists. Given $\Gamma$ uniformizing $(0;p,p,p,p)$, we have a tiling
$\Delta(\alpha_0,\beta_0,\gamma_0)$ as before. Let $R_\pi := \pi^{-1}(A) \cup \pi^{-1}(B) \cup \pi^{-1}(C) \cup \pi^{-1}(D)$. We say that
$R_\pi$ is composed by four disjoint classes. Let $d(X,Y)$ be the hyperbolic distance in $\Delta$ between the points $X$ and $Y$. Consider a
(hyperbolic) triangle from this tiling, i.e., having interior angles $\alpha_0$, $\beta_0$ and $\gamma_0$. If all of them are $\leq
\frac{\pi}{p}$, then there is nothing to prove. Otherwise, suppose that $\gamma_0> \frac{\pi}{p}$. Fix $A' \in \pi^{-1}(A)$ and consider the
quadrilateral with vertices $\{A',B',D'',C' \}$ in Figure \ref{f2}. Now, we join by an hyperbolic segment the points $C'$ and $D'$ and erase the
hyperbolic segment which joins $A'$ and $B'$. Do this throughout $\Delta(\alpha_0,\beta_0,\gamma_0)$ and obtain a new tiling for the same
$\Gamma$ given by $\Delta(\alpha_1,\beta_1,\gamma_1)$ ($\alpha_1 + \beta_1 + \gamma_1 = \frac{2 \pi}{p}$). By applying the hyperbolic sine
theorem, we notice that one of the edges of the new triangle is strictly smaller than the corresponding edge of the old one. This new smaller
edge has one vertex equal to $A'$. Now we ask whether we satisfy the condition of the proposition or not, but notice that it could not be
satisfied only by the angles whose vertex is not $A'$. If one of them is $>\frac{\pi}{p}$, we repeat this process. Hence, since we are getting a
chain of inequalities $ ... < d(A',P_n) < ... < d(A',P_1) < d(A',P_0) $ where $P_i \in \pi^{-1}(B) \cup \pi^{-1}(C) \cup \pi^{-1}(D)$ and
$R_\pi$ is a discrete set, this process must stop. For the uniqueness statement, the key part is to show that the smallest distances between the
four classes of $R_\pi$ are achieved by the edges of the triangle having interior angles $\alpha$, $\beta$ and $\gamma$ with $\alpha, \beta,
\gamma \leq \frac{\pi}{p}$. This can be worked out by using basic hyperbolic geometry on $\Delta$.

\begin{obs}
We can also define a triangular tiling on $\C$ to study compact Riemann surfaces of genus one (this is done in section \ref{s8}). In that case,
we can prove the same proposition and even better we have that two surfaces are isomorphic if and only if the corresponding tilings having acute
angles are the same. When the genus is greater than two, we get a combinatorial problem to solve. \label{r1}
\end{obs}

From now on, we choose for $\Gamma$ the tiling $\Delta(\alpha,\beta,\gamma)$ as in Prop. \ref{p6}; if one of the angles is $\frac{\pi}{p}$, then
we take $\Delta(\frac{\pi}{p},\beta,\gamma)$ such that $\beta \geq \gamma$. With this, we have assigned a unique $\Delta(\alpha,\beta,\gamma)$
to each $\Gamma$, this will be the $\bm{canonical \ tiling}$ for $\Gamma$. Now, we notice that Prop. \ref{p2} also applies to the case of
surfaces in $\gimel_p$. The following theorem is a consequence of Prop. \ref{p2} plus Prop. \ref{p6}.

\begin{teo}
Let $S,S' \in{\gimel_p}$ and suppose that $S \simeq S'$. Then, they have the same canonical tiling. \label{t3}
\end{teo}

\section{uniformization}\label{s5}

We fix $\Gamma$ and the corresponding canonical tiling $\Delta(\alpha,\beta,\gamma)$. We recall that in section $3$ we chose a fundamental
domain $\rm{P}(\alpha,\beta,\gamma)$ for $\Gamma^*$. To obtain the surface $S = \Delta/ \Gamma^*$ we will need to specify the identifications of
the edges of this polygon $\rm{P}(\alpha,\beta,\gamma)$. This polygon is formed by $p$ hexagons, identical to $\rm{H}_0$, around $A'$. We
denoted those hexagons by $\rm{H}_0, \rm{H}_1,...,\rm{H}_{p-1}$. Let us write $b_1 := \overline{D'B'}$, $b_2 := \overline{B'D''}$, $c_1 :=
\overline{D''C'}$ and $c_2 := \overline{C'D'''}$; similarly for $\rm{H}_{k}$ we have $b_{2k+1}$, $b_{2k+2}$, $c_{2k+1}$ and $c_{2k+2}$ with $k
\in{ \{0,1,...,p-1 \}}$, obtained by rotating $\rm{H}_0$ counterclockwise $k$ times and letting those edges be the images of $b_1$, $b_2$, $c_1$
and $c_2$ respectively. We notice that by symmetry, we only need to say how to glue $b_1$ and $c_1$ with some other edges. When we glue $w_i$
with $w_j$, we will write $w_i \leftrightarrow w_j$.

\begin{prop}
The only possibilities for edge identification are $b_1\leftrightarrow b_{2i+2}$, $c_1\leftrightarrow c_{2j+2}$ with $i,j \in{\{1,2,...,p-1\}}$,
excluding when $i+j+1\equiv 0 (mod \ p)$. \label{p7}
\end{prop}

Considering all the possible identifications is equivalent to considering all the possibilities for surjective group homomorphisms $\phi :
\Gamma \longrightarrow \langle u : u^p=1 \rangle$ such that $ker(\phi)=\Gamma^*$ and $\phi(x)=u$ (we are using the notation of section $3$).
These possibilities are given by $\phi(y)=u^i$ and $\phi(z)=u^j$ such that $i,j \in{\{1,...,p-1\}}$ and $p$ does not divide $i+j+1$ (the last
condition is equivalent to $\phi(xyz) \neq 1$).

The last proposition gave the rules for identifications, but does not say anything about whether different identifications will give
analytically distinct surfaces. Actually, we saw in section \ref{s2} that for some distinct identifications we can obtain isomorphic Riemann
surfaces. That problem will be solved later. For now, we want to compute $\Gamma^*$ for a given $b_1\leftrightarrow b_{2i+2}$ and
$c_1\leftrightarrow c_{2j+2}$. The automorphism of $\Gamma^*$ which identifies an edge $w_n$ with an edge $w_m$ will be denoted by $T_{w_n
\leftrightarrow w_m}$. By solving linear equations module $p$, which come from the constraints $x^i=y$ and $x^j=z$ in $\Gamma/\Gamma^*$, we can
prove the following proposition.

\begin{prop}
If $b_1\leftrightarrow b_{2i+2}$ and $c_1\leftrightarrow c_{2j+2}$, then $\Gamma^*$ is generated by $\{x^{ni}T_{b_1\leftrightarrow b_{2i+2}}
x^{-ni},x^{nj}T_{c_1\leftrightarrow c_{2j+2}} x^{-nj} \}_{n=0}^{p-1}$, where $T_{b_1\leftrightarrow b_{2i+2}} = x^i y^{-1}$ and
$T_{c_1\leftrightarrow c_{2j+2}} = x^j z^{-1}$. \label{p8}
\end{prop}

\section{special domains and the main theorem}\label{s6}

As we did with Lefschetz surfaces, we now want to consider some special fundamental domains for $\Gamma^*$. The fundamental domain
$P(\alpha,\beta,\gamma)$ given in section \ref{s3} will be denoted by $[AD]$ since it is centered at $A'$ ($\pi(A')=A$) and has $2p$ vertices in
$\pi^{-1}(D)$. Similar to the case of Lefschetz surfaces, let us denote by $[XY]$ a domain hyperbolically identical to $P(\alpha,\beta,\gamma)$
but now centered at $X'$ such that $\pi(X')=X$ and has $2p$ vertices in $\pi^{-1}(Y)$, where $X,Y \in{ \{A,B,C,D \} }$. Hence, we have only
twelve possible special domains: $[AD]$, $[AB]$, $[AC]$, $[BC]$, $[BA]$, $[BD]$, $[CA]$, $[CB]$, $[CD]$, $[DA]$, $[DB]$ and $[DC]$. Proposition
\ref{p8} gives us the gluing for $[AD]$ encoded by a pair of numbers $(i,j)$. Since we also want to keep track of the geometry of the domain
$[XY]$, we will include a subindex which indicates the angle at the vertex $Y''$ opposite to $X'$ (e.g., for $[AD]$ in Figure \ref{f2}, this
angle is $\alpha$). We denote the angles by numbers: $\alpha:=1$, $\beta:=2$ and $\gamma:=3$. In this way, for $[AD]$ we have $(i,j)_1$. For any
special domain $[XY]$, we denote this pair by $s[XY]$ (we will refer to it as a pair, even though it has a third number given by the angle). The
following proposition gives the pairs for all the special domains.

\begin{prop}
Suppose that $s[AD]=(i,j)_1$ and $k=-(i+j+1)$, then: $s[AB]=(j,k)_3$, $s[AC]=(k,i)_2$, $s[BD]=(i^{-1}j,i^{-1})_2$, $s[BC]=(i^{-1},i^{-1}k)_1$,
$s[BA]=(i^{-1}k,i^{-1}j)_3$, $s[CD]=(j^{-1},ij^{-1})_3$, $s[CA]=(ij^{-1},j^{-1}k)_2$, $s[CB]=(j^{-1}k,j^{-1})_1$, $s[DA]=(jk^{-1},ik^{-1})_1$,
$s[DC]=(ik^{-1},k^{-1})_3$ and $s[DB]=(k^{-1},jk^{-1})_2$, where all the numbers involved are taken from $\Z/p\Z$ in $\{ 1,2,...,p-1 \}$.
\label{p9}
\end{prop}

Now to each surface $S \in \gimel_p$ with a fixed $\Delta(\alpha,\beta,\gamma)$, we can assign the following set:
\bigskip

- When the canonical tiling is neither $\Delta(\frac{2\pi}{3p},\frac{2\pi}{3p},\frac{2\pi}{3p})$ nor
$\Delta(\frac{\pi}{p},\frac{\pi}{2p},\frac{\pi}{2p})$, \\
$ \Lambda_{(i,j)}^p := \{(i,j)_1,(j,k)_3,(k,i)_2,(i^{-1}j,i^{-1})_2,(i^{-1},i^{-1}k)_1,(i^{-1}k,i^{-1}j)_3,$\\
$(j^{-1},ij^{-1})_3,(ij^{-1},j^{-1}k)_2,(j^{-1}k,j^{-1})_1,(jk^{-1},ik^{-1})_1,(ik^{-1},k^{-1})_3,(k^{-1},jk^{-1})_2\}$.\\

- (\bm{$Equilateral \ tiling$}) If $\Delta(\frac{2\pi}{3p},\frac{2\pi}{3p},\frac{2\pi}{3p})$, \\
$ ^{e}\Lambda_{(i,j)}^p :=
\{(i,j), (j,k),(k,i),(i^{-1}j,i^{-1}),(i^{-1},i^{-1}k),(i^{-1}k,i^{-1}j),(j^{-1},ij^{-1}),$\\
$(ij^{-1},j^{-1}k),(j^{-1}k,j^{-1}),(jk^{-1},ik^{-1}),(ik^{-1},k^{-1}),(k^{-1},jk^{-1})
\}$.\\

- (\bm{$Square \ tiling$}) If $\Delta(\frac{\pi}{p},\frac{\pi}{2p},\frac{\pi}{2p})$, $^c\Lambda_{(i,j)}^p := \Lambda_{(i,j)}^p \cup
\Lambda_{(j,i)}^p$.

\begin{teo}
$S, S' \in{\gimel_p}$ are isomorphic if and only if $\Delta(\alpha,\beta,\gamma)=\Delta(\alpha',\beta',\gamma')$ and
$\Lambda_{(i,j)}^p=\Lambda_{(i',j')}^p$. \label{t4}
\end{teo}

The proof is the proof of Thm. \ref{t1} adapted to our situation.

\section{the numbers}\label{s7}
We first want to study the set $^e\Lambda_{(i,j)}^p$. For that, we define
$$ \Sigma_p := \{ (i,j)\in{\Z/ p\Z \times \Z/ p\Z} \ \ : \ \ i\neq 0 \ , \ j\neq 0 \ , \ i+j+1\neq 0 \} \ .$$
One can check that the sets $^e\Lambda_{(i,j)}^p$ partition $\Sigma_p$.

\begin{prop}
Let $p>3$. The following are all the possible cases for $^e\Lambda_{(i,j)}^p$,
\begin{itemize}
\item[-] $^e\Lambda_{(1,1)}^p= \{ (1,1),(1,-3),(-3,1),(-3^{-1},-3^{-1})\}$. \item[-] $^e\Lambda_{(1,p-1)}^p= \{ (1,-1),(-1,-1),(-1,1) \}$.
\item[-] $^e\Lambda_{(i,p-1)}^p= \{(i,-1),(1,-i),(-i,i),(-i^{-1},i^{-1}),(i^{-1},-1),(-1,-i^{-1})\}$ with $i\neq \pm 1$. \item[-] $12$ different
pairs otherwise.
\end{itemize}
For $p=3$, we have only one case $^e\Lambda_{(2,2)}^p= \{ (2,2),(2,1),(1,2) \} = \Sigma_3$. \label{p10}
\end{prop}

We remark that two sets $^e\Lambda_{(i,j)}^p$ are either equal or disjoint, since they partition $\Sigma_p$. This gives an equivalence relation,
and the class of $^e\Lambda_{(i,j)}^p$ will be denoted again by $^e\Lambda_{(i,j)}^p$. When we add the subindex to study $\Lambda_{(i,j)}^p$, we
need to consider the set $\Sigma_p$ with the subindices $\{ 1,2,3 \}$, and this gives a partition of the corresponding new set. Again, we will
denote the class of $\Lambda_{(i,j)}^p$ in this new equivalence relation by $\Lambda_{(i,j)}^p$. The classes of this new partition will be
divided in cases $\bm{\kappa_l}$ in the following way:

\begin{itemize}
\item $\kappa_1$: Here we only have $\Lambda_{(1,1)}^p$. \item $\kappa_2$: This is formed by the classes $\Lambda_{(1,i)}^p$,
$\Lambda_{(i,1)}^p$ with $i \in{\{2,3,...,p-4\}}$, and $\Lambda_{(i,i)}^p$ with $i \in{\{2,3,...,p-2\}}$ and $i\neq (p-3)^{-1}$. \item
$\kappa_3$: Here we only have the classes of $\Lambda_{(-1,-1)}^p$. \item $\kappa_4$: This is formed by the classes of $\Lambda_{(-1,i)}^p$ with
$i\in{\{2,3,...,p-2\}}$ and $p \nmid i^2+1$. \item $\kappa_5$: This only happens when $p\equiv 1(mod \ 4)$ and the classes are
$\Lambda_{(-1,i)}^p$ with $i^2\equiv -1(mod \ p)$. \item $\kappa_6$: All the rest.
\end{itemize}

In the appendix, we give tables for the sets $\Lambda_{(i,j)}^p$ for some primes, divided according to these cases. By counting the different
classes in $\Sigma_p$, we obtain the following proposition.

\begin{prop}
Let $p>3$. Then, the number of analytically distinct Riemann surfaces in $\gimel_p$ with equilateral tiling
$\Delta(\frac{2\pi}{3p},\frac{2\pi}{3p},\frac{2\pi}{3p})$ is $\frac{p^2+11}{12}$. If $p=3$, then there is only one. \label{p11}
\end{prop}

\begin{prop}
The number of analytically distinct Riemann surfaces in $\gimel_p$ with a fixed tiling which is neither equilateral
$\Delta(\frac{2\pi}{3p},\frac{2\pi}{3p},\frac{2\pi}{3p})$ nor square $\Delta(\frac{\pi}{p},\frac{\pi}{2p},\frac{\pi}{2p})$ is $\frac{p^2+3}{4}$.
\label{p12}
\end{prop}

\begin{proof}
For this calculation, we notice that given $(i,j)_n$, the pairs which have the same subindex are $(i^{-1},i^{-1}k)_n$, $(j^{-1}k,j^{-1})_n$ and
$(jk^{-1},ik^{-1})_n $, and so we have three different classes which contain pairs of this type. A priori, we could think that we have three
different surfaces for each class. However, these classes can intersect themselves. This intersection can only happen in the case $\kappa_1$.
Then, by using the previous calculation in Prop. \ref{p11}, we obtain $ \frac{12}{4} \times (\frac{p^2+11}{12}-1) + 1 = \frac{p^2+3}{4} $
surfaces if $p$ is not $3$. For $p=3$, we get $3$ non-isomorphic surfaces.
\end{proof}

\begin{prop}
The number of analytically distinct Riemann surfaces in $\gimel_p$ with a square tiling $\Delta(\frac{\pi}{p},\frac{\pi}{2p},\frac{\pi}{2p})$ is
$\frac{p^2+2p+5}{8}$ if $p\equiv 1(mod \ 4)$, and $\frac{p^2+2p+1}{8}$ if $p\equiv -1(mod \ 4)$. \label{p13}
\end{prop}

The proof is similar to the proof of Prop. \ref{p12}, with the difference that we need to include more classes since the set
$^c\Lambda_{(i,j)}^p$ is formed by more elements.

\section{a parameterization space for $\gimel_p$}\label{s8}

The goal now is to build a parameterization space for $\gimel_p$. This will be a suitable quotient of $\frac{p^2+3}{4}$ copies of $\C$ at the
points representing square and equilateral tilings. In this parameterization space, two distinct points will be two distinct surfaces of
$\gimel_p$. The set $\gimel_2$ is the set of all tori, and a parameterization for it is the well known moduli space of curves of genus one $M_1
\cong \C$. Consider $M_1$ as the quotient of the upper half plane $\H$ by the group $PSL(2,\Z)$. This group has as a fundamental domain the
region $\mathcal{R}$ given by $x\geq0$, $x\leq1$, $x^2+y^2\geq1$ and $(x-1)^2+y^2\geq1$ where $x$ and $y$ are the real coordinates for points in
$\H$. Each point $\tau$ in $\mathcal{R}$ gives a unique triangle formed by $\tau$ and $1$ having acute interior angles. If we take the triangles
given by $r i$ and $(r+1)i$ ($r \in \R$), then they differ only by reordering the angles. If we consider $\tau$ and $\tau'$ in $x^2+y^2 = 1$ and
$(x-1)^2+y^2 = 1$ such that there is $g \in PSL(2,\Z)$ with $g(\tau)=\tau'$, then the corresponding isosceles acute triangles are the same. In
conclusion, we have a natural bijection between acute triangles and tori (up to the case when we have a right triangle, in which case we take
triangle($\frac{\pi}{2},\alpha,\beta)=$ triangle($\frac{\pi}{2},\beta,\alpha$)). This was exactly what we did with canonical tilings. So, we
naturally obtain a bijection between $\C$ and the set of canonical tilings. In this way, we think about the space of canonical tilings as $\C$.
Actually, this is true because canonical tilings are in one to one correspondence with the moduli space of four unordered points in $\P^1$ (up
to projective equivalence), which is isomorphic to $\C$.

Now, to build a parameterization space for $\gimel_p$, we take $\frac{p^2+3}{4}$ copies of $\C$ and quotient at the points where the tiling is
equilateral or square according to the cases $\kappa_l$. We will get three topologically different types of connected components. This is done
below (we take $p>3$).

\begin{center}
\begin{figure}
\includegraphics[width=8.5cm]{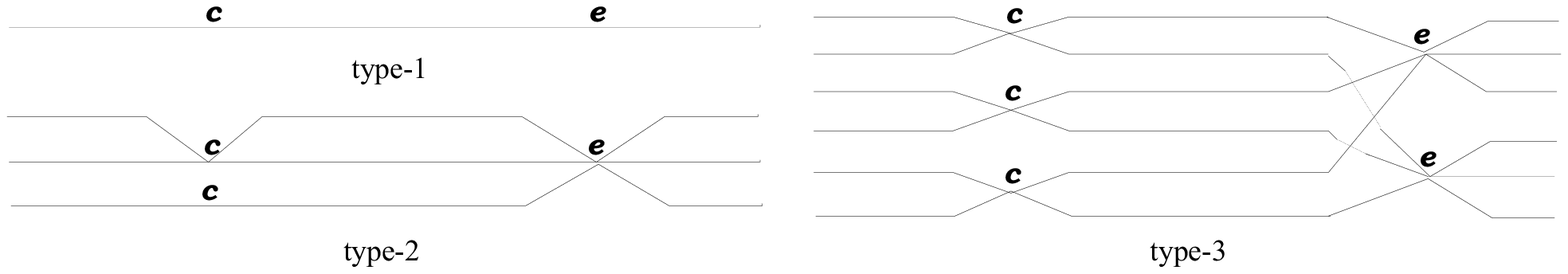}
\caption{Component types.} \label{f3}
\end{figure}
\end{center}

\begin{center}
\begin{figure}
\includegraphics[width=9cm]{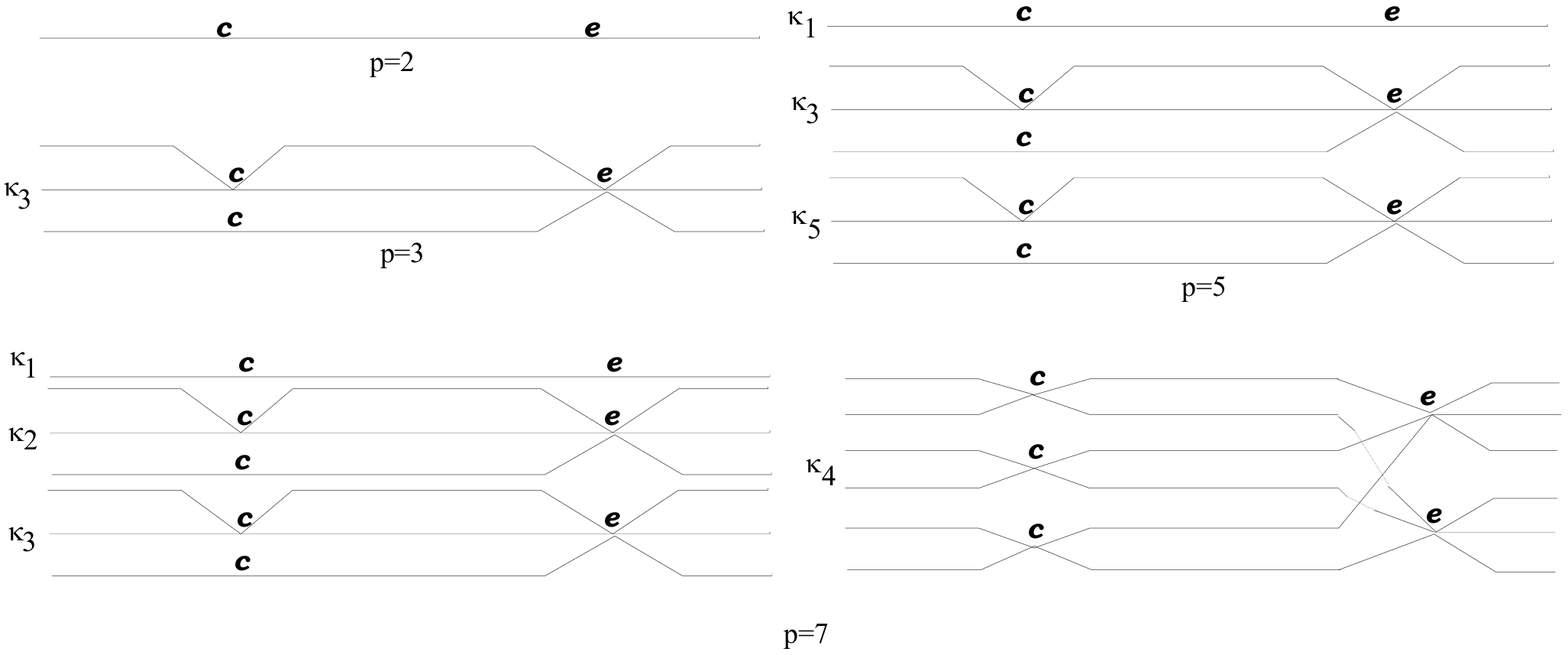}
\caption{$\gimel_2$, $\gimel_3$, $\gimel_5$ and $\gimel_7$. $e$ means equilateral tiling and $c$ means square tiling.} \label{f4}
\end{figure}
\end{center}

$\kappa_1$: We have only one $\C$ for every $p$ (\textsf{component type-1}).

$\kappa_2$ and $\kappa_3$ and $\kappa_5$: Each of these cases will be formed by $3$ copies of $\C$ gluing two in the square point and the third
at the equilateral point, as in Figure \ref{f3} (\textsf{component type-2}). In total, for these cases we have $\frac{p-3}{2}$ disjoint
components if $p\equiv -1(mod \ 4)$, and $\frac{p-1}{2}$ otherwise.

$\kappa_4$ and $\kappa_6$: Each component is formed by $6$ copies of $\C$, gluing three of them with the other three at the square point (so we
have three square surfaces) and, by taking one from each square pair, we glue each three of them at the equilateral point, as shown in Figure
\ref{f3} (\textsf{component type-3}). In total, for these cases we have $\frac{p^2-6p+5}{24}$ disjoint components if $p\equiv 1(mod \ 4)$, and
$\frac{p^2-6p+17}{24}$ otherwise.

For $p=3$, we only have one component type-2. In this way, for $p>3$, we have described $\gimel_p$ as a topological space divided in several
connected components. The number of connected components is: one if $p=3$, $\frac{p^2+6p+17}{24}$ if $p\equiv 1(mod \ 4)$, and
$\frac{(p+5)(p+1)}{24}$ if $p\equiv -1(mod \ 4)$. In Figure \ref{f4}, we show the parametrization spaces of $\gimel_p$ for some primes .

\section{full automorphism groups}\label{s9}

Let $p$ be an odd prime. In this section we are going to compute the full automorphism group for each surface in $\gimel_p$. Let $S$ be in
$\gimel_p$ and let $T_S \in{\aut(S)}$ be an element of order $p$ (prime as always) fixing four points. We define $\aut(S)'=\{ f \in{\aut(S)} \ :
\ f(\fix(T_S))=\fix(T_S) \ \}$. Then, it can be checked that this subgroup of $\aut(S)$ is the normalizer of $\langle T_S \rangle$ in $\aut(S)$.
What we will do is to compute $\aut(S)'$ and then say when $\langle T \rangle$ is actually normal in $\aut(S)$. We are able to compute
$\aut(S)'$ by using the sets $\Lambda_{(i,j)}^p$ as vectors in the sense that we keep track of the special domains together with the pairs of
numbers. The idea is that we have an extra automorphism (apart from $T_S$) if and only if there are two special domains with the same pair of
numbers (i.e., same edge identification plus the same angle). Then, we look at all the possible combinations among the special domains and
relations between the new automorphisms, so that we finally compute this group. For instance, if all the domains have different pairs of
numbers, then $\aut(S)'$ will be $\Z/p\Z$, i.e., $S$ will not have any extra automorphism fixing the set $\fix(T_S)$. All the proofs go in that
way. Below we show a list of propositions which give $\aut(S)'$ for every $S \in{\gimel_p}$, according to the cases $\kappa_l$ and the canonical
tilings. We denote by $\bm{^c\gimel_p}$ the set of surfaces with square tiling and $\bm{^e\gimel_p}$ the set of surfaces with equilateral
tiling.

\begin{prop}
Let $S$ be in $\kappa_1$, then: $S \in{^e\gimel_p}$ implies $\aut(S)' \simeq \Z/3p \Z$, $S \in{^c\gimel_p}$ implies $\aut(S)' \simeq \Z/2p \Z$
and for all the rest we have $\aut(S)' \simeq \Z/p \Z$. A surface $S$ in $\kappa_1$ is never hyperelliptic. \label{p15}
\end{prop}

\begin{prop}
Let $S$ be in $\kappa_2$, then: \\
- If $S \in{^e\gimel_p}$, $\aut(S)' \simeq \Z/ p\Z$ .\\
- If $S \in{^c\gimel_p}$ and $S$ belongs to the class of $\Lambda_{(i,i)}^p$, then $\aut(S)' \simeq \Z/ 2p\Z$, otherwise $\aut(S)' \simeq \Z/ p\Z$. \\
- Otherwise, $\aut(S)' \simeq \Z/ p\Z$.\\
There are not hyperelliptic surfaces in this case. \label{p16}
\end{prop}

\begin{prop} (hyperelliptic component of $\gimel_p$) Let $S$ be in $\kappa_3$, then:\\
- If $S \in{^e\gimel_p}$, $\aut(S)' \simeq D_{2p}$.\\
- If $S \in{^c\gimel_p}$ and $S$ belongs to the class of $\Lambda_{(-1,-1)}^p$, then $\aut(S)' \simeq (\Z/ 2\Z \times \Z/2p\Z) \rtimes \Z/ 2\Z$,
otherwise $\aut(S)' \simeq D_{2p}$. \\
- For all the rest, $\aut(S)' \simeq D_{2p}$.\\
In this case, all the surfaces are hyperelliptic. \label{p17}
\end{prop}

\begin{prop}
Let $S$ be in $\kappa_4$, then $\aut(S)' \simeq D_{p}$. There are not hyperelliptic surfaces in this case. \label{p18}
\end{prop}

\begin{prop}
Let $S$ be in $\kappa_5$ (only when $p\equiv 1(mod \ 4)$), then:\\
- If $S \in{^e\gimel_p}$, $\aut(S)' \simeq D_{p}$.\\
- If $S \in{^c\gimel_p}$ and $S$ belongs to the class of $\Lambda_{(i,-i)}^p$, then $\aut(S)' \simeq \Z/ p\Z \rtimes \Z/ 4\Z$, otherwise
$\aut(S)' \simeq D_{p}$.\\
- For all the rest, $\aut(S)' \simeq D_{p}$.\\
There are not hyperelliptic surfaces in this case. \label{p19}
\end{prop}

\begin{prop}
Let $S$ be in $\kappa_6$, then $\aut(S)' \simeq \Z/ p\Z$. \label{p20}
\end{prop}

In \cite{S}, Singerman worked out a list for finitely maximal Fuchsian groups. This has been used several times to compute automorphism groups
for compact Riemann surfaces (for example see \cite{B-C} or \cite{K-S}). By using those methods, it can be proved that $\langle T_S \rangle$ is
a normal subgroup of $\aut(S)$ for $p>5$. The proof can be found in \cite{M}, \cite{H-M} or \cite{W}. Now, for $p=3$ or $p=5$, we have explicit
lists with the classification of the automorphisms groups. For $p=3$ it is classical, and for $p=5$ (so genus $4$), it can be found in
\cite{K-K}. Hence, we can check that the only cases when $\langle T_S \rangle$ is not normal are:\\
- The unique surface of genus $2$ which has $48$ automorphisms ($\aut \simeq GL(2,\F_3)$) and is given in $\C^2$ by $y^2=x(x^4-1)$.\\
- The unique surface of genus $4$ which has $120$ automorphisms ($\aut \simeq S_5$) (permutation group of five elements). This is the famous
Bring's
curve (see \cite{Ri-R-2}).\\

Therefore, up to the two cases above, the previous list of propositions classify the full automorphism groups for all the surfaces in
$\gimel_p$.

\section{Components of the singular locus of $M_g$}\label{s10}

In this section, we compute the number of components of dimension zero and one of the singular locus $\bm{S_g}$ of $M_g$, the moduli space of
smooth curves of genus $g$. The dimension zero case was computed in \cite{G}. It is known (for example see \cite{O}) that for $g\geq 4$, $S_g =
\{ [C] \in M_g : \vert\aut(C)\vert \neq 1 \} $. In what follows, we will use the notation given in \cite{C} and \cite{G}. The next result
appears in those papers.

\begin{teo}
Let $p$ be a prime number, $\{ a_1,..., a_n \}$ natural numbers with $a_i<p$ and $\sum_{i=1}^n a_i \equiv 0(mod \ p)$; $g'\leq g$ natural
numbers with $2g-2=p(2g'-2)+n(p-1)$. There exists a closed subscheme $S(p,g';a_1,...,a_n)$ of the singular locus $S_g$ irreducible of dimension
$3g'-3 +n$, whose geometric points are the curves of genus $g$ that are coverings of the form $\overline{C}(p,X'; \sum_{i=1}^n a_i q_i, L)$ for
some curve $X'$ of genus $g'$, $\{q_1,...,q_n\}$ different points of $X'$ and $L$ an invertible sheaf on $X'$ such that $L^p \simeq
\mathcal{O}_{X'}(\sum_{i=1}^n a_1 q_i)$. \label{t5}
\end{teo}

The translation to our situation is the following. The number $n$ is the amount of fixed points of the action of $\Z/p\Z$ on a compact Riemann
surface $C$ and $g'=0$. Hence, for different elections of the numbers $a_l$'s, $S(p,0;a_1,a_2,a_3)$ are the Lefschetz surfaces and $S(p,0;
a_1,a_2,a_3,a_4)$ are the $\gimel_p$ surfaces. For $n=3$, we can choose to have $a_1=1$, $a_2=k$ and $a_3=p-1-k$; and for $n=4$, we can choose
to have $a_1=1$, $a_2=i$, $a_3=j$ and $a_4=k$ with the notation of section \ref{s6}. For different elections of the tuples of numbers $a_l$'s,
we can have the same subscheme $S(p,0;a_1,...,a_n)$. In \cite{G-G}, G. Gonz\'alez defines equivalence classes for those tuples, denoting them by
$\overline{m}$. He proves that $M_g^p$, which is the union of all the sets $S(p,0;a_1,...,a_n)$ in $M_g$ for a fixed $n$, is a disjoint union of
normal irreducible subvarieties $M_g^p(\overline{m})$, where $M_g^p(\overline{m})$ is $S(p,0;a_1,...,a_n)$ in the equivalent class given by
$\overline{m}$. He states that the number of components of $M_g^p$ can be read from the generating function given by Lloyd in \cite{Ll}. On the
other hand, we computed the components for Lefschetz surfaces (which is just the number of them for a fixed $p$) and for $\gimel_p$ (see last
part of section \ref{s8}). We compared our result with the corresponding number from the generic formula in \cite{Ll}, and we got the same.

Now, by Thm. \ref{t5}, we can say that the dimension zero and the dimension one components of $S_g$ are contained in the Lefschetz surfaces and
$\gimel_p$ respectively. To find them, we need to subtract the components which have generically the action of an extra automorphism, other than
the one which produces the action of $\Z/p\Z$. In another words, we need to subtract the components of Riemann surfaces having $\aut \neq
\Z/p\Z$. But we already know by sections \ref{s2} and \ref{s9} which they are. Putting all together, we obtain the following theorems.

\begin{teo}
When $g\geq 4$, the number of isolated singularities of $M_g$ is $\frac{g-2}{3}$ if $2g+1$ is a prime number and $2g+1\equiv 2(mod \ 3)$, or
$\frac{g-3}{3}$ if $2g+1$ is a prime number and $2g+1\equiv 1(mod \ 3)$. If $g=2$ or $3$, then $M_g$ has only one isolated singularity.
\label{t6}
\end{teo}

\begin{teo}
The number of dimension one components of the singular locus of $M_g$ is $\frac{g(g+2)}{24}$ if $g+1>3$ is a prime number, or zero otherwise.
\label{t6}
\end{teo}

We observe that Thm. \ref{t6} does not agree with the tables given by Cornalba in \cite{C} for dimension one. For example, in his tables he has
$S(7,0;1,2,5,6)$ as a component for $S_6$, but we proved in section \ref{s9} that $S(7,0;1,2,5,6)$ has generically $D_7$ acting, and actually it
is contained in $S(2,3;1,1)$, which is a component of $S_6$ of dimension $8$.  What we think is that he included in the components our cases
$\kappa_4$ and $\kappa_5$, but we proved in section \ref{s9} that generically they have $D_p$ acting, and so those cases are contained in bigger
dimensional varieties of $S_g$. If we take out those cases, we coincide with the numbers in his tables.\footnote{Cornalba agreed with this
observation and he will publish an erratum to "On the locus of curves with automorphisms" (\cite{C}) in the corresponding journal.}

We would like to notice that our classification also says how these components intersect $S_g$. This is contained in section \ref{s9}. For
example, the component for the $\kappa_1$ case intersects $S_g$ only at two points: the corresponding square Riemann surface with $\aut \simeq
\Z/2p\Z$ and the corresponding equilateral Riemann surface with $\aut \simeq \Z/3p\Z$. We can compute all the intermediate coverings given by
the extra elements of order $2$ or $3$ respectively. Also, as a curiosity, we can imply that $M_{p-1}$ has a unique curve with $\aut \simeq
\Z/3p\Z$. Same thing can be done with all the dimension  zero and one components.

\section{some special families in $\gimel_p$ and their affine equations}\label{s11}

In the present section we will give equations for all Riemann surfaces which correspond to special tilings, i.e., equilateral tilings
$\Delta(\frac{2\pi}{3p},\frac{2\pi}{3p},\frac{2\pi}{3p})$ and square tilings $\Delta(\frac{\pi}{p},\frac{\pi}{2p},\frac{\pi}{2p})$, and also for
Hyperelliptic Riemann surfaces.
\bigskip

{\emph Hyperelliptic surfaces in $\gimel_p$}

In the previous sections, we found that for each $\gimel_p$ we have a connected component of type-2 (see section \ref{s8}) where every surface
is hyperelliptic. Outside of this component, we do not have hyperelliptic surfaces. This component corresponds to the case $\kappa_3$
represented by (and suitable permutations of the subindices)
$$ \Lambda_{(-1,-1)}^p = \{(-1,-1)_1,(-1,1)_3,(1,-1)_2,(1,-1)_2,(-1,-1)_1,$$
$$ (-1,1)_3, (-1,1)_3,(1,-1)_2,(-1,-1)_1,(-1,-1)_1,(-1,1)_3,(1,-1)_2 \} \ . $$

\begin{prop}
Let $S\in{\gimel_p}$ be hyperelliptic. Then, $S$ can be represented by the affine curve in $\C^2$, $ y^2 = (x^p-a^p)(x^p+\frac{1}{a^p})$ where
$a \in{\C}$ is nonzero. \label{p21}
\end{prop}

For instance, this can be proved by using the above calculation for $\aut(S)'$ which indicates how the ramification points of the hyperelliptic
involution can be moved by an element in $\aut(S)$. According to the equation above, we can calculate the generators of $\aut(S)$, which are
$(x,y) \mapsto (-\frac{1}{x},\frac{iy}{x^p})$ and $(x,y) \mapsto (e^{\frac{2\pi i}{p}}x,-y)$. For $a=\pm1$, we will add to the list $(x,y)
\mapsto (-x,y)$. Section \ref{s8} tells us the following.

\begin{prop}
Let $p>2$ and $S\in{\gimel_p}$ be hyperelliptic. Then, excepting for the curve of $48$ automorphisms in $\gimel_3$, we have:\\
- If $S$ is the unique Riemann surface in $^c\Lambda_{(-1,-1)}^p$ $(y^2=x^{2p}-1)$, then $$\aut(S) \simeq (\Z/ 2\Z \times \Z/2p\Z) \rtimes \Z/
2\Z$$ - Otherwise, $\aut(S) \simeq D_{2p}$. \label{p22}
\end{prop}

\begin{obs}
Our work implies that for each $p>3$ prime number, we have a unique Riemann surface $S$ in $M_{p-1}$ such that $\aut(S) \simeq (\Z/ 2\Z \times
\Z/ 2p\Z) \rtimes \Z/ 2\Z$ (this is a big group according to the usual definition, see \cite{A}). This generalizes the unique curve of genus $2$
having $24$ automorphisms. \label{o1}
\end{obs}

\begin{obs}
Let $N(g)$ be the maximal order of the full automorphism group of curves in $M_g$. It is known (see \cite{A}) that $8(g-1) \leq N(g) \leq
84(g-1)$. Our work tell us that if $\vert \aut(S) \vert = N(g)$ and $g\geq 5$, then either $S$ is the unique curve with $\aut(S) \simeq (\Z/ 2\Z
\times \Z/ 2(g-1)\Z) \rtimes \Z/ 2\Z$ or for every prime divisor $p$ of $N(g)$ we have $p\leq g$. For genus two, three and four, it is known
that $N(g)$ is only satisfied by the curve of $48$ automorphisms, the Klein curve and the Bring curve respectively. We notice that those curves
are the only cases where $p>g$ and $\Z/p\Z$ is not normal in their full automorphism group. \label{o2}
\end{obs}

{\emph{Equilateral surfaces in $\gimel_p$}

- For each $p>3$, we have $\frac{p^2+11}{12}$ non isomorphic equilateral surfaces, and for $p=3$ only one.

- With respect to $\aut(S)$, we have:
\begin{itemize}
\item[a)] $\aut(S_{(1,1)}) \simeq \Z/ 3p\Z$. $S_{(1,1)}$ is the unique surface with this group for each $p>3$. \item[b)] $\aut(S_{(-1,-1)})
\simeq D_{2p}$ (one for each $p>2$). \item[c)] $\aut(S) \simeq D_{p}$ (there are $\frac{p-3}{2}$ for each $p$). \item[d)] $\aut(S) \simeq \Z/
p\Z$ (all the rest).
\end{itemize}

- The affine equations as curves in $\C^2$ are $ y^p =(x-1)(x-w)^n(x-w^2)^m $ where $w= e^{\frac{2 \pi i}{3}}$ ($i=\sqrt{-1}$),
$n,m\in{\{1,2,...,p-1}\}$ and $p\nmid (n+m+1)$. This is because of our rotation numbers and the fact that there is a unique curve in $\gimel_p$
with $\Z/3\Z$ acting. The classification according to $\aut(S)$ is:

\begin{itemize}
\item[a)] $ y^p =(x-1)(x-w)(x-w^2)= x^3-1$ for $\aut(S_{(1,1)}) \simeq \Z/ 3p\Z$.

\item[b)] $ y^p =(x-1)(x^2+x+1)^{p-1}$ (Hyperelliptic).

\item[c)] $y^p = (x-1)(x-w)^{n}(x-w^2)^{p-1}$, $n\in{\{2,...,p-2\}}$, $\aut(S_{(n,-1)}) \simeq D_{p}$.

\item[d)] All the rest.
\end{itemize}

All the generators for these groups can be computed using section \ref{s8} and the corresponding equation (this also holds for the square case
below). For instance, $\aut(S_{(1,1)}) = \langle (x,y) \mapsto (e^{\frac{2 \pi i}{3}}x, e^{\frac{2 \pi i}{p}}y) \rangle \simeq \Z/3p\Z$.
\bigskip

{\emph{Square surfaces in $\gimel_p$}

- Let $p$ be odd. Then, the number of square surfaces is $\frac{p^2+2p+5}{8}$ if $p\equiv 1(mod \ 4)$, or $\frac{p^2+2p+1}{8}$ if $p\equiv
-1(mod \ 4)$.

- Below, we give a list of the square surfaces in $\gimel_p$ according to the $\kappa$ case. If $S=S_{(a,b)_1} \in{^c\gimel_p}$, then $S$ is
represented by the equation in $\C^2$: $y^p=(x-1)(x-i)^a(x+1)^c(x+i)^b$ with $c=2p-1-a-b$. Again, we know this because of our interpretation of
the rotation numbers and the existence of a unique surface with full automorphism group isomorphic to $(\Z/ 2\Z \times \Z/ 2p\Z) \rtimes \Z/
2\Z$.

\begin{itemize}
\item[$\kappa_1$:] This case has only one surface for each $p>3$ and $\aut(S_{(1,1)_1}) \simeq \Z/ 2p\Z$ with equation
$y^p=(x-1)(x^2+1)(x+1)^{p-3}$.

\item[$\kappa_2$:] In this case we have several surfaces and two possibilities for $\aut(S)$, depending on whether $(a,b)_1$ has $a=b$ or not.

\begin{itemize}
\item[1)] If $S_{(a,a)_1}$, then $\aut(S_{(a,a)_1}) \simeq \Z/2p\Z$ and its equation is $ y^p=(x-1)(x^2+1)^a(x+1)^{2p-1-2a}$. \item[2)] In the
other case, we have the equation evaluated in the respective $(a,b)_1$ and $\aut(S) \simeq \Z/ p\Z$.
\end{itemize}

\item[$\kappa_3$:] The hyperelliptic case gives us two possible $\aut(S)$, with the respective equations $y^p=(x^2-1)(x^2+1)^{p-1}$ for $(\Z/2\Z
\times \Z/ 2p\Z) \rtimes \Z/ 2\Z$ and $y^p=(x-1)(x+i)(x-i)^{p-1}(x+1)^{p-1}$ for $D_{2p}$.

\item[$\kappa_4$:] In this case, we have only one possible group: $\aut(S) \simeq D_p$ and the equation is evaluated in the respective numbers.

\item[$\kappa_5$:] This case is restricted to $p\equiv 1(mod \ 4)$ and we have two possibilities:

\begin{itemize}
\item[1)] If $S_{(a,-a)_1}$ with $a^2\equiv -1(mod \ p)$, then $\aut(S_{(a,-a)_1}) \simeq \Z/ p\Z \rtimes \Z/ 4\Z$ and its equation is
$y^p=(x-1)(x-i)^a(x+1)^{p-1}(x+i)^{p-a}$.

\item[2)] It is analogous to the case $\kappa_4$.
\end{itemize}

\item[$\kappa_6$:] The rest.
\end{itemize}

\begin{obs}
Our work implies that for each $p>3$ prime number such that $p\equiv 1(mod \ 4)$, we have a unique Riemann surface $S$ in $M_{p-1}$ with
$\aut(S) \simeq \Z/ p\Z \rtimes \Z/ 4\Z$. This surface generalizes the Bring's curve (\cite{Ri-R-2}) in the sense that this curve is $S$ for
$p=5$. Hence, it is given by the equation $y^5= (x-1)(x-i)^2(x+1)^3(x+i)^4 $. We remark that for $p>5$, the corresponding curve in $\gimel_p$
has $\Z/p\Z$ normal in $\aut(S)$, and this is not the case for the Bring's curve which has full automorphism group isomorphic to $S_5$.
\label{o1}
\end{obs}
\bigskip

\section{appendix}\label{s12}

\begin{center}
\textbf{$\Omega_k^p$ for $p \in{ \{5,7,11,13,17,19 \} }$}
\end{center}

$$ \Omega_1^5=\{1,2,3\} $$

$$ \Omega_1^7=\{1,3,5\} \  \Omega_2^7=\{2,4\} $$

$$ \Omega_1^{11}=\{1,5,9\} \  \Omega_2^{11}=\{2,3,4,6,7,8\} $$

$$ \Omega_1^{13}=\{1,6,11\} \  \Omega_2^{13}=\{2,10,5,7,8,4\} \  \Omega_3^{13}=\{3,9\} $$

$$ \Omega_1^{17}=\{1,8,15\} \  \Omega_2^{17}=\{2,5,7,9,11,14\} \   \Omega_3^{17}=\{3,4,6,10,12,13\} $$

$$ \Omega_1^{19}=\{1,9,17\} \ \Omega_2^{19}=\{2,6,8,10,12,16\}  \  \Omega_3^{19}=\{3,4,5,13,14,15\}  \ \Omega_7^{19}=\{7,11\}$$

\newpage

\begin{center}
\textbf{$\Lambda_{(i,j)}^p$ for $p \in{ \{5,7,11,13 \} }$}
\end{center}
\begin{center}
\tiny{
\begin{tabular}{|r|c|r|r|r|c|r|r|r|c|r|r|r|c|r|r|r|c|r}
\hline
            & p=3 &     &     & p=5 &     &     &     &     & p=7 &     &     \\
            &     &     &     &     &     &     &     &     &     &     &     \\ \hline
   $[AD]_1$ & 2.2 &     & 1.1 & 1.4 & 2.3 &     & 1.1 & 1.2 & 1.6 & 2.5 & 3.4 \\ \hline
   $[AB]_3$ & 2.1 &     & 1.2 & 4.4 & 3.4 &     & 1.4 & 2.3 & 6.6 & 5.6 & 4.6 \\ \hline
   $[AC]_2$ & 1.2 &     & 2.1 & 4.1 & 4.2 &     & 4.1 & 3.1 & 6.1 & 6.2 & 6.3 \\ \hline
            &     &     &     &     &     &     &     &     &     &     &     \\ \hline
   $[BD]_2$ & 1.2 &     & 1.1 & 4.1 & 4.3 &     & 1.1 & 2.1 & 6.1 & 6.4 & 6.5 \\ \hline
   $[BC]_1$ & 2.2 &     & 1.2 & 1.4 & 3.2 &     & 1.4 & 1.3 & 1.6 & 4.3 & 5.2 \\ \hline
   $[BA]_3$ & 2.1 &     & 2.1 & 4.4 & 2.4 &     & 4.1 & 3.2 & 6.6 & 3.6 & 2.6 \\ \hline
            &     &     &     &     &     &     &     &     &     &     &     \\ \hline
   $[CD]_3$ & 2.1 &     & 1.1 & 4.4 & 2.4 &     & 1.1 & 4.4 & 6.6 & 3.6 & 2.6 \\ \hline
   $[CA]_2$ & 1.2 &     & 1.2 & 4.1 & 4.3 &     & 1.4 & 4.5 & 6.1 & 6.4 & 6.5 \\ \hline
   $[CB]_1$ & 2.2 &     & 2.1 & 1.4 & 3.2 &     & 4.1 & 5.4 & 1.6 & 4.3 & 5.2 \\ \hline
            &     &     &     &     &     &     &     &     &     &     &     \\ \hline
   $[DA]_1$ & 2.2 &     & 3.3 & 1.4 & 2.3 &     & 2.2 & 3.5 & 1.6 & 2.5 & 3.4 \\ \hline
   $[DC]_3$ & 2.1 &     & 3.3 & 4.4 & 3.4 &     & 2.2 & 5.5 & 6.6 & 5.6 & 4.6 \\ \hline
   $[DB]_2$ & 1.2 &     & 3.3 & 4.1 & 4.2 &     & 2.2 & 5.3 & 6.1 & 6.2 & 6.3 \\ \hline
            &     &     &     &     &     &     &     &     &     &     &     \\
            &$\kappa_3$&     &$\kappa_1$&$\kappa_3$&$\kappa_5$&
            &$\kappa_1$&$\kappa_2$&$\kappa_3$&$\kappa_4$&$\kappa_4$     \\ \hline
\end{tabular}
}
\end{center}

\begin{center}
\tiny{
\begin{tabular}{|r|c|r|r|r|c|r|r|r|c|r|r|r|c|r|r|r|c|r}
\hline
            &     &     &     &     &     & p=11 &     &     &     &     &      \\
            &     &     &     &     &     &      &     &     &     &     &      \\ \hline
   $[AD]_1$ & 1.1 & 1.2 & 1.3 & 1.4 & 1.10& 2.3  & 2.5 & 2.10& 8.10& 3.10& 5.10 \\ \hline
   $[AB]_3$ & 1.8 & 2.7 & 3.6 & 4.5 &10.10& 3.5  & 5.3 & 10.9& 10.3& 10.8& 10.6 \\ \hline
   $[AC]_2$ & 8.1 & 7.1 & 6.1 & 5.1 &10.10& 5.2  & 3.2 & 9.2 & 3.8 & 8.3 & 6.5  \\ \hline
            &     &     &     &     &     &      &     &     &     &     &      \\ \hline
   $[BD]_2$ & 1.1 & 2.1 & 3.1 & 4.1 &10.10& 7.6  & 8.6 & 5.6 & 4.7 & 7.4 & 2.9  \\ \hline
   $[BC]_1$ & 1.8 & 1.7 & 1.6 & 1.5 & 1.10& 6.8  & 6.7 & 6.10& 7.10& 4.10& 9.10 \\ \hline
   $[BA]_3$ & 8.1 & 7.2 & 6.3 & 5.4 &10.10& 8.7  & 7.8 & 10.5& 10.4& 10.7& 10.2 \\ \hline
            &     &     &     &     &     &      &     &     &     &     &      \\ \hline
   $[CD]_3$ & 1.1 & 6.6 & 4.4 & 3.3 &10.10& 4.8  & 9.7 & 10.9& 10.3& 10.8& 10.6 \\ \hline
   $[CA]_2$ & 1.8 & 6.9 & 4.2 & 3.4 &10.10& 8.9  & 7.5 & 9.2 & 3.8 & 8.3 & 6.5  \\ \hline
   $[CB]_1$ & 8.1 & 9.6 & 2.4 & 4.3 & 1.10& 9.4  & 5.9 & 2.10& 8.10& 3.10& 5.10 \\ \hline
            &     &     &     &     &     &      &     &     &     &     &      \\ \hline
   $[DA]_1$ & 7.7 & 5.8 & 6.2 & 3.9 & 1.10& 5.7  & 9.8 & 6.10& 7.10& 4.10& 9.10 \\ \hline
   $[DC]_3$ & 7.7 & 8.8 & 2.2 & 9.9 &10.10& 7.9  & 8.4 & 10.5& 10.4& 10.7& 10.2 \\ \hline
   $[DB]_2$ & 7.7 & 8.5 & 2.6 & 9.3 &10.10& 9.5  & 4.9 & 5.6 & 4.7 & 7.4 & 2.9  \\ \hline
            &     &     &     &     &     &      &     &     &     &     &     \\
            &$\kappa_1$&$\kappa_2$&$\kappa_2$&$\kappa_2$&$\kappa_3$&$\kappa_6$
            &$\kappa_6$&$\kappa_4$&$\kappa_4$&$\kappa_4$&$\kappa_4$     \\ \hline
\end{tabular}
}
\end{center}

\begin{center}
\tiny{
\begin{tabular}{|r|c|r|r|r|c|r|r|r|c|r|r|r|c|r|r|r|c|r}
\hline
          &     &     &     &     &     &      &     & p=13&     &     &     &     &     &     &      \\
          &     &     &     &     &     &      &     &     &     &     &     &     &     &     &      \\ \hline
 $[AD]_1$ & 1.1 & 1.2 & 1.3 & 1.4 & 1.5 & 1.12 & 2.3 & 2.4 & 2.11& 2.12& 3.4 & 3.5 & 3.10& 3.12& 5.8  \\ \hline
 $[AB]_3$ & 1.10& 2.9 & 3.8 & 4.7 & 5.6 & 12.12& 3.7 & 4.6 &11.12&12.11& 4.5 & 5.4 &10.12&12.10& 8.12 \\ \hline
 $[AC]_2$ & 10.1& 9.1 & 8.1 & 7.1 & 6.1 & 12.1 & 7.2 & 6.2 & 12.2& 11.2& 5.3 & 4.3 & 12.3& 10.3& 12.5 \\ \hline
          &     &     &     &     &     &      &     &     &     &     &     &     &     &     &      \\ \hline
 $[BD]_2$ & 1.1 & 2.1 & 3.1 & 4.1 & 5.1 & 12.1 & 8.7 & 2.7 & 12.7& 6.7 & 10.9& 6.9 & 12.9& 4.9 & 12.8 \\ \hline
 $[BC]_1$ & 1.10& 1.9 & 1.8 & 1.7 & 1.6 & 1.12 & 7.10& 7.3 & 7.6 & 7.12& 9.6 & 9.10& 9.4 & 9.12& 8.5  \\ \hline
 $[BA]_3$ & 10.1& 9.2 & 8.3 & 7.4 & 6.5 & 12.12& 10.8& 3.2 & 6.12& 12.6& 6.10& 10.6& 4.12& 12.4& 5.12 \\ \hline
          &     &     &     &     &     &      &     &     &     &     &     &     &     &     &      \\ \hline
 $[CD]_3$ & 1.1 & 7.7 & 9.9 &10.10& 8.8 & 12.12& 9.5 & 10.7& 6.12&12.11& 10.4& 8.11& 4.12&12.10& 5.12 \\ \hline
 $[CA]_2$ & 1.10& 7.11& 9.7 & 10.5& 8.9 & 12.1 & 5.11& 7.8 & 12.7& 11.2& 4.11& 11.6& 12.9& 10.3& 12.8 \\ \hline
 $[CB]_1$ & 10.1& 11.7& 7.9 & 5.10& 9.8 & 1.12 & 11.9& 8.10& 7.6 & 2.12&11.10& 6.8 & 9.4 & 3.12& 8.5  \\ \hline
          &     &     &     &     &     &      &     &     &     &     &     &     &     &     &      \\ \hline
 $[DA]_1$ & 4.4 & 6.3 & 2.5 & 8.2 & 3.11& 1.12 & 6.4 & 5.9 & 2.11& 7.12& 6.11& 11.4& 3.10& 9.12& 5.8  \\ \hline
 $[DC]_3$ & 4.4 & 3.3 & 5.5 & 2.2 &11.11& 12.12& 4.2 & 9.11&11.12& 12.6& 11.8& 4.10&10.12& 12.4& 8.12 \\ \hline
 $[DB]_2$ & 4.4 & 3.6 & 5.2 & 2.8 & 11.3& 12.1 & 2.6 & 11.5& 12.2& 6.7 & 8.6 &10.11& 12.3& 4.9 & 12.5 \\ \hline
          &     &     &     &     &     &      &     &     &     &     &     &     &     &     &      \\
          &$\kappa_1$&$\kappa_2$&$\kappa_2$&$\kappa_2$&$\kappa_2$&$\kappa_3$
          &$\kappa_6$&$\kappa_6$&$\kappa_4$&$\kappa_4$&$\kappa_6$&$\kappa_6$
          &$\kappa_4$&$\kappa_4$&$\kappa_5$     \\ \hline
\end{tabular}
}
\end{center}
\bigskip


\bigskip
Department of Mathematics, University of Michigan. \\ 530 Church Street Ann Arbor, MI 48109-1043.

\end{document}